\documentclass[12pt]{article}
\newcommand{\nocopyright}{
No Copyright\thanks{
The authors hereby waive all copyright
and related or neighboring rights to this work,
and dedicate it to the public domain.
This applies worldwide.
}}
\title{Maybe there's no such thing as a random sequence}
\author{
Peter G. Doyle\thanks{Dartmouth College.}
}
\date{Version 1.0 dated 17 March 2010
\\ \nocopyright
}
\usepackage{amsthm,amssymb}
\usepackage{verbatim}
\usepackage{amsmath}
\usepackage{graphicx}

\newtheorem{definition}{Definition}
\newtheorem{theorem}{Theorem}

\newcommand{\ZF}{\mathrm{ZF}}
\newcommand{\ZFC}{\mathrm{ZFC}}
\newcommand{\SM}{\mathrm{SM}}
\newcommand{\Con}{\mathrm{Con}}
\begin{document}
\maketitle

\begin{abstract}
An infinite binary sequence is deemed to be random
if it has all definable properties that hold almost surely
for the usual probability measure on the set of infinite binary sequences.
There are only countably many such properties, so
it would seem that the set of random sequences
should have full measure.
But in fact there might be no random sequences,
because for all we know, there might be no undefinable sets.
\end{abstract}

\centerline{\emph{For Laurie and Jim}}

\section{What is a random sequence?}

About 30 years ago now, my friend and mentor J. Laurie Snell got
interested in the question of what constitutes an infinite random sequences
of $0$s and $1$s.
The model here is an infinite sequence of independent flips of a fair
coin,
with heads counting as $1$ and tails as $0$.
Or more technically, the standard product measure on
$\prod_{i=1}^\infty \{0,1\}$,
which up to a little fussing is the same as the standard Borel
measure on the unit interval $[0,1]$.

Some sequences are obviously not random:
\begin{itemize}
\item
$0000000000\ldots$;
\item
$1111111111\ldots$;
\item
$0101010101\ldots$;
\item
$1101001000100001\ldots$.
\end{itemize}
Other presumably non-random sequences are the binary expansion
of $1/\pi$, or $1/e$, or any sequence that can be printed out
by a computer program.
As there are only a countable number of computer programs,
this gives us only a countable set;
the complement is still an uncountable set of full measure.

Throwing out these computable sequences doesn't come near to doing the job.
Many more sequences need to be weeded out.
For example, a random sequence should exhibit the
strong law of large numbers:  Asymptotically it should have half $0$s and half
$1$s.
There are uncountably many sequences that fail this test,
e.g. all those of the form $00a00b00c00d00e\ldots$
(On Beyond Zebra!).

And having half $0$s and half $1$s is far from enough.
Really we want the sequence to have all those properties dictated by
probability theory,
like say,
the law of the iterated logarithm.
This is a relative of the central limit theorem
which states---well, you can look it up.
What matters is that it is a statement about the sequence that
either is true or false,
and which according to probability theory is true \emph{almost surely}.

Now, some people have felt that a test for randomness should be effective in
some sense.
For example, maybe you should be able to make money from a non-random sequence,
as you could for example from a sequence of more than half $1$s by betting
on $1$ each time.
Looking into this led Laurie into a thicket of papers by Kolmogorov,
Martin-L\"of,
Schnorr, and Chaitin,
which are all very interesting and answer the question of what is random
from a certain perspective.
(This approach is by now a thriving industry:
See Downey et al.\ \cite{downeyetal:calibrating};
Nies
\cite{nies:randomness};
and
Downey and Hirschfeldt
\cite{downeyhirschfeldt:algorithmic}.)

At this point Jim Baumgartner, a logician who had been drawn into this
morass by Laurie, proposed the following:

\begin{definition} \label{def:random}
A sequence is \emph{non-random}
if and only if it belongs to some definable set of sequences
having measure $0$.
Here definable means uniquely definable by a formula in the language of first order set theory having only one free variable (no parameters).
And by measure $0$ we mean outer measure $0$,
in case the set is not Borel-measurable.
\end{definition}

Since there are only a countable number of formulas,
we're throwing out  a countable collection of sets of measure $0$,
so the sequences that remain---the \emph{random} sequences,
should form a set of full measure.
Almost every sequence should be random.

I accepted this definition of random sequence for over 30 years.
Then about a month ago,
thanks to a stimulating colloquium talk by Johanna Franklin
and subsequent discussions with Rebecca Weber and Marcia Groszek,
I came to realize that
the argument that there are plenty of random sequences is bogus.
(Or perhaps it would be better to say that it is `suspect', in light
the fact that the same reasoning was used by Tarski---see section
\ref{sec:bogus} below.)
It is possible that
under this definition
\emph{there are no random sequences}!
The reason is that
the standard of axioms of set theory
(assuming they are consistent)
do not rule out the possibility that \emph{every set is definable}.
There are models of set theory, satisfying all the standard axioms,
where every set in the universe is definable.
(Cf.\ section \ref{sec:pdm} below.)
Of course these are countable models.
In these models, every sequence is definable,
and hence every singleton set consisting of a single sequence.
Since a singleton set has measure $0$,
every sequence belongs to a definable set of measure $0$, hence is
non-random.
So in such a model there are no random sequences.

Now we might have objected earlier,
`Of course there are no random sequences:
Given any sequence $\sigma$, the singleton set
$\{ \sigma \}$ has measure $0$.
A random sequence cannot be equal to any particular sequence.'
We thought we were avoiding this by only considering
definable sets of sequences:
There are only a countable number of definable sequences,
so there should be a full measure set of sequences left after we've ruled out
definable sequences.
And even after we've gone on to throw out
all definable sets of sequences of measure $0$,
there should remain a full-measure set of random sequences.

But now we're saying that in fact this objection might be
justified, after all:
For all we know, it might be the case that \emph{all sequences are definable}.
So maybe there really are no random sequences.

\section{How can this be?}

The answer is \emph{Skolem's paradox}.
We're dealing here with a countable model of set theory.
In any model of set theory, the collection of definable sets
will be countable from outside the model.
If the model is countable,
there is no obvious impediment to having every set be definable.
And in fact this turns out to be possible.

Look, the real problem here is that you can't define definability.
There is no formula in the language of set theory
characterizing a definable set,
because there is no formula characterizing a true formula.
(If there were, we'd be in real trouble.)
And so we can't write a formula characterizing
random sequences.
What we can do is write a formula characterizing random sequences
\emph{within a given model of set theory}.
Now, if the standard axioms of set theory are consistent,
then there is a model of set theory satisfying these axioms,
and in this model every set could be definable, which would
mean in particular that there are no random sequences in the model.

Note that we are not saying that it must be the case that there are
no random sequences.
There certainly are models where not all sequences are definable
(assuming set theory is consistent).
There presumably are models where the random sequences have full measure.
We're just saying that it \emph{may be}
that there are no random sequences.

\section{What to make of all this?}

One sensible response would be that we have missed the boat.
This definition of random sequence
is too restrictive.
We should be less demanding.
We should climb into the boat with Kolmogorov,
Martin-L\"of, Schnorr, and Chaitin.
And then we can talk not just about infinite sequences, but finite sequences
as well.
No finite sequences are completely random, of course, but clearly some are more
random than others.

We prefer to stick with Definition \ref{def:random},
and consider the possibility that there really are no random sequences.
Maybe the Old Man doesn't play dice with the universe of sets.

That is, assuming there really is a universe of sets.
In this connection, I can't resist quoting Abraham Robinson
\cite[p.\ 230]{robinson:64}:
\begin{quotation}
My position concerning the foundations of Mathematics is based
on the following two main points or principles.

(i) Infinite totalities do not exist in any sense of the word
(i.e., either really or ideally).
More precisely, any mention,
or purported mention,
of infinite totalities is, literally,
\emph{meaningless}.

(ii) Nevertheless, we should continue the business of Mathematics
``as usual,'' i.e., we should act \emph{as if}
infinite totalities really existed.
\end{quotation}
See also Cohen \cite{cohen:comments}.

\section{Models of set theory for which all sets are definable} \label{sec:pdm}

I've found that Cohen's book `Set Theory and the Continuum Hypothesis'
\cite{cohen:continuum} is
a good place to look for general background on model theory.
The book is addressed to non-specialists, and `emphasizes
the intuitive motivations while at the same time giving as complete
proofs as possible'.

Here,
quoted verbatim from Cohen
\cite[pp.\ 104--105]{cohen:continuum},
are precise statements about models where all sets are definable.

\begin{theorem}
$\ZF + \SM$ implies the existence of a unique transitive model $M$
such that if $N$ is any standard model there is an $\in$-isomorphism
of $M$ into $N$.  $M$ is countable.
\end{theorem}

Here $\SM$ is the statement that $\ZF$ has a standard model,
meaning one where the membership relation in the model coincides with
the `real world' membership relation $\in$.
The existence of a standard model is not provable because it implies
$\Con(\ZF)$.
Cohen
\cite[p.\ 79]{cohen:continuum} says that $\SM$ is `most probably ``true"',
and gives an intuitive argument for accepting it as an axiom.
$\SM$ holds just if it is possible to quit early in the transfinite induction
that produces Godel's constructible universe $L$, and still have a model of
$\ZF$.
To get the minimal model $M$, we stop the construction at the earliest possible
ordinal.
$M$ satisfies $V=L$
(cf. \cite[p.\ 104]{cohen:continuum}),
hence also the axiom of choice, so in $M$ we have
a standard model of $\ZFC + (V=L)$.

\begin{theorem}
For every element $X$ in $M$ there is a formula $A(y)$ in $\ZF$ such
that $x$ is the unique element in $M$ satisfying $A_M(x)$.
Thus in $M$ every element can be `named'.
\end{theorem}

Models where every set is definable are called `pointwise definable'.
So according to these results, if $\ZF$ has a standard model, then
$\ZFC$ has a pointwise definable standard model.

The requirement that $\ZF$ has a standard model can be dispensed with,
if we don't care about winding up with a pointwise definable model
that is non-standard.
John Steel points out that
the definable sets within any model $N$ of $\ZFC + (V=L)$
constitute an elementary submodel $H$.
This implies that $H$ is a model of $\ZFC + (V=L)$,
and every set in $H$ is definable in $H$ (not just in $N$).
So, starting with any model for $\ZF$, we can restrict to a model
of $\ZFC + (V=L)$,
and within that find a pointwise definable model of $\ZFC +(V=L)$.

For much more about pointwise definable models of set theory, see
Hamkins, Linetsky, and Reitz \cite{hamkinsetal:pointwise}.

\section{What was wrong with the proof?} \label{sec:bogus}

So, what was wrong with the proof that there are plenty of random sequences?
Let's look at the proof, given by Tarski
\cite[p.\ 220]{tarski:definissable}
in 1931,
that there exist undefinable sets of real numbers.
This same proof method would show that there exist undefinable real
numbers,
or what is the same, undefinable random sequences.
Tarski's original paper
\cite{tarski:definissable} is in French.
Here,
from Tarski
\cite[p.\ 119]{tarski:definable},
is an English translation:

\begin{quote}
Moreover it is not difficult to show that the family of all definable
sets (as well as that of the functions which determine them) is only
denumerable, while the family of \emph{all} sets of numbers is not
denumerable.  The existence of undefinable sets follows immediately.
\end{quote}

If anyone but Tarski had written this, I think we would say that
the author is confusing the system with the metasystem.  In the
metasystem, we can talk about the family $D$ of definable sets of reals
of the system, and prove that it is denumerable \emph{in the metasystem}.
In the system itself, we can't talk about $D$, but we can talk about the
family $P$ of all sets of reals of the system, and prove that it is not
denumerable \emph{within the system itself}.
But from this we cannot
conclude that $D$ differs from $P$:
$D$ is denumerable in the metasystem; 
$P$ is not denumerable in the system.
There is no contradiction here.

Tarski goes on to give a second proof:
\begin{quote}
`Plus encore' [the translation has `Also', but a better rendering 
might be `And not only that'],
the definable sets can be arranged in an ordinary infinite sequence;
by applying the diagonal procedure it is possible \emph{to define
in the metasystem a concrete set which would not be definable in the
system itself}.  In this there is clearly no trace of any antinomy,
and this fact will not appear at all paradoxical if we take proper
note of the relative character of the notion of definability.
\end{quote}

If anyone but Tarski had written this, I think we would say that the
author has failed to take proper note of the relative character of
\emph{being a set}.
This concrete set that we've defined in the metasystem
might not be a set in the system itself.  So we have failed to produce
an undefinable set of the system.

Now, this critique of Tarski's reasoning has been written from the perspective
of model theory---a modern perspective based in large part
on work that Tarski did after writing down these proofs.
We interpret the `metasystem'
as a formal system in which statements about models (`systems')
of set theory can be formulated and proven.
If we choose as the metasystem ZFC, the usual formal system for set
theory, then in this metasystem we can prove a statement to the effect
that there exist models of ZFC where every set is definable.
This means that Tarski's arguments that there must be undefinable sets
cannot be correctly formalized in ZFC (assuming ZFC is consistent).
And we have pointed out where the problem lies in trying to formalize them.

So, does this mean that Tarski's arguments are \emph{wrong}?
They are if it is fair to recast them in the framework of model theory.
But it is not really clear that this is fair.
It has been suggested, for example, that Tarski is thinking of the
`standard model' rather than some arbitrary model.
That's all very well, but what does it mean concretely?
What special methods of proof apply to the `standard model'?

If we can make sense of Tarski's arguments,
then presumably we can salvage the argument
that the set of random sequences from Definition \ref{def:random}
has full measure.
Maybe it is false that maybe there are no random sequences.

\section*{Acknowledgement and disclaimer}
I'm grateful for help I've received from
Jim Baumgartner, Johanna Franklin,
Steven Givant, Marcia Groszek, Joel David Hamkins,
Laurie Snell, John Steel, and Rebecca Weber.
None of them is responsible for my opinions, or whatever errors
may be on display here.
Despite the impression I've been trying to give,
I'm aware that I know next to nothing about logic.
I've written this note because I find this whole business intriguing,
and I believe other people will too.

\bibliography{random}
\bibliographystyle{plain}
\end{document}